\documentclass{amsart}
\usepackage{amssymb}

\newtheorem{theorem}{Theorem}[section]

\newtheorem{lemma}[theorem]{Lemma}
\newtheorem{proposition}[theorem]{Proposition}
\newtheorem{remark}[theorem]{Remark}

\begin{document}

\title{Group gradings on finite dimensional Lie algebras}
\author[D.~Pagon, D.~Repov\v s, and M.~Zaicev]{Du\v san Pagon, Du\v san Repov\v s, and Mikhail Zaicev}

\address{Du\v{s}an~Pagon \\
Faculty of Natural Sciences and Mathematics, University of
Maribor, Gosposvetska 84, Maribor, 2000, Slovenia}
\email{dusan.pagon@uni-mb.si}

\address{Du\v san Repov\v s \\Faculty of Mathematics and Physics, and
Faculty of Education, University of Ljubljana,
P.~O.~B. 2964, Ljubljana, 1001, Slovenia}
\email{dusan.repovs@guest.arnes.si}

\address{Mikhail Zaicev \\Department of Algebra\\ Faculty of Mathematics and
Mechanics\\  Moscow State University \\ Moscow,119992 Russia\\
zaicevmv@mail.ru}

\thanks{This research was supported by the Slovenian Research Agency grants
P1-0292-0101, J1-2057-0101, and J1-9643-0101.
The third author was partially supported by RFBR grants No
 09-01-00303, 09-0190416-Ukr-f-a  and SSC-1983.2008.1}



\begin{abstract}
We study gradings by noncommutative groups on finite dimensional Lie
algebras over an algebraically closed field of characteristic zero. 
It is shown that if $L$ is gradeg by a non-abelian finite group $G$
then the solvable radical $R$ of $L$ is $G$-graded and there exists
a Levi subalgebra $B=H_1\oplus\cdots\oplus H_m$ homogeneous in $G$-grading
with graded simple summands $H_1, \ldots, H_m$. All supports $Supp~H_i,
i=1\ldots, m$, are commutative subsets of $G$.
\end{abstract}

\maketitle

\section{Introduction}

Graded rings and graded algebras 
have been
extensively
studied  during the last decades (cf. e.g.
\cite{BP2009} - \cite{K1968}).
Group gradings were investigated both in the
associative case (\cite{BZ2002} - 
\cite{CM1984}, \cite{S1997}) and in the Lie case
(\cite{BS2002},\cite{K1968},\cite{K1972},\cite{PZ89},\cite{S1999}),
or in other nonassociative cases (\cite{BP2009},\cite{BS2002},\cite{E1998}).

One of the important tasks is the description of all possible gradings on different algebras.
For example, one of the well-known and actively used results in Lie theory is the
description of $\mathbb{Z}$-gradings on finite-dimensional complex Lie algebras \cite{K1972}.
The description of all $\mathbb{Z}_2$-gradings on matrix algebras plays an
exceptional role in the theory of algebras with polynomial identities
(cf.  \cite{GZ2005}).

All abelian gradings on matrix algebras were described in \cite{BSZ2001}.
It was found that the so-called fine gradings on full matrix algebra
play an exceptional role in the theory of orthogonal Cartan
decompositions of simple complex Lie algebras
(cf. \cite{KT1992}).
For the noncommutative case all group gradings on
matrices were described in \cite{BZ2002}. In particular,
it was shown that in general case fine gradings are closely
connected with some problems of the theory of projective representations of finite groups.

 All finite dimenional graded simple associative algebras were described in \cite{BZS2008}. This description is presently 
 being
 very intensively applied
 in PI-theory (cf. e.g. \cite{AL}, \cite{KZ}).
All finite dimensional $\mathbb{Z}_2\times \mathbb{Z}_2$-graded simple Lie algebras
were described in \cite{BP2009}
and this result was applied for the classification of simple color Lie superlagebras.
So, the
description
of
all possible gradings on algebras plays an important role both 
in
the
structure theory of finite dimensional and infinite dimensional graded algebras and 
in
its applications.

In the present paper paper we shall
study general properties of commutative and noncommutative gradings
on finite dimensional Lie algebras over an algebraically closed field
of characteristic zero. First, we shall
show that the description of noncommutative gradings on semisimple algebras
can be reduced to abelian gradings and to the classification of graded simple algebras
(cf. Propositions \ref{p1} and \ref{p2}.)
Then we shall characterize finite dimensional graded simple algebras and
show that all of them can be split
into four classical series $A,B,C,D$ and five exceptional series,
according to the
classification of simple Lie algebras (cf.
Proposition \ref{p2} and Remark \ref{r1}). Finally, we shall
show that in the
case of finite groups,
any graded Lie algebra is a split extension
of a
homogeneous Levi subalgebra and
a
solvable radical (cf. Theorem \ref{t1}).

\section{General properties of noncommutative gradings}

Let $G$ be a group. Given a Lie algebra $L$ over a field $F$, we say that $L$ is $G$-graded if it can be decomposed into a
direct sum of subspaces
$$
L=\bigoplus_{g\in G} L_g
$$
such that $[L_g, L_h]\subseteq L_{gh},$ for any $g,h\in G$.
Subspaces $L_g$ are called homogeneous components and an element $x\in L$
is called homogeneous if $x\in L_g,$ for some $g\in G$.
In this case we write for shortness $\deg x=g$.
A subspace (resp. 
subalgebra, ideal) $V$ is said to be a graded subspace (resp.
subalgebra, ideal) if
$$
V=\bigoplus_{g\in G} (V\cap L_g).
$$
In other words,
if $x=x_{g_1}+\cdots+x_{g_n},$
where $\deg x_{g_i}=g_i$ and $g_1,\ldots, g_n\in G$
are pairwise distinct,
then $x\in V$ if and only if all $x_{g_1},\ldots,x_{g_n}$ belong to $V$.

A subset
$$
Supp~L=\{g\in G\vert L_g\ne 0\}
$$
is called the support of the
grading. Actually,
we can suppose that $Supp~L$ generates $G$.
We call an algebra $L$ graded simple if $[L,L]\ne 0$ and $L$ does not
contain
any
non-trivial graded ideals.

We shall
usually use the notation $[x_1,\ldots,x_n] $
for 
the
left-normed product in
a
Lie algebra.
That is,
$[x_1,\ldots,x_n]=[[x_1,\ldots,x_{n-1}],x_n]$ for all $n\ge 3$. Similarly,
we
shall
write $[W_1,\ldots,W_n]=[W_1,\ldots,W_{n-1},W_n]$ for any subspaces $W_1,\ldots, W_n\subseteq L$.

It is well-known that in the
case of a
simple Lie algebra $L=\oplus_{g\in G}L_g,$
the group $G$ must be commutative
(cf. e.g.
\cite{PZ89}).
First, we generalize this property to the
case of graded simple Lie algebras. In the following
two lemmas and Proposition \ref{p1},
$G$ is an arbitrary group
and $L$ is not necessarily
a
finite dimensional Lie algebra over an arbitrary field $F$.

\begin{lemma}\label{l1}
Let $L=\oplus_{g\in G}$ be a $G$-graded Lie algebra and
let $[L_{g_1},\ldots,L_{g_m}] \ne 0,$
for some $g_1,\ldots,g_m\in G$. Then $g_1,\ldots,g_m$ commute in $G$.
\end{lemma}

{\em Proof}. First, note that the inequality $[L_g,L_h]\ne 0$ implies $gh=hg$ since $[L_g,L_h]=[L_h,L_g]\subseteq L_{gh}\cap L_{hg}$.
Hence our statement is obvious for $m=2$.
Now we apply the induction on $m$. Suppose that
$[x_1,\ldots, x_m]\ne 0,$ for some $x_1\in L_{g_1}, \ldots, x_m\in L_{g_m}$ and $m\ge 3$.
Since $[x_1,\ldots, x_m]=[u, x_{m-1}] + [v,[x_{m-1}, x_m]],$ where $u=[x_1,\ldots,x_{m-2}, x_m]$ and $v=[x_1,\ldots, x_{m-2}]$ (or $v=x_1$ in case $m=3$)
either $[u, x_{m-1}]\ne 0$ or $[v,[x_{m-1}, x_m]]\ne 0$.

In the first case $g_1,\ldots, g_{m-2},g_m$ commute by the inductive hypothesis
and $g_{m-1}$ commutes with the product $g_1\cdots g_{m-2}g_m$. Since $g_1,\ldots,g_{m-1}$
also commute, by  induction we have $g_mg_{m-1}=g_{m-1}g_m$ that is,
$g_m$
commutes with all $g_1,\ldots,g_{m-1}$.

In the second case, $g_1,\ldots,g_{m-1}$ commute and  the product $g_{m-1}g_m$ commutes
with all $g_i, 1\le i\le m-2$. Since $g_ig_{m-1}=g_{m-1}g_i$, $1\le i\le m-2$, we obtain
$g_ig_{m}=g_{m}g_i,$ for all $i=1,\ldots, m-2$. Clearly, $g_mg_{m-1}=g_{m-1}g_m$ and we
have thus
completed
the proof.
\hfill $\Box$

\begin{lemma}\label{l2}
Let $L=\oplus_{g\in G}L_g$ be a $G$-graded Lie algebra and
let $gh\ne hg$ for some $g,h\in Supp~L$. Then $I=[Id(L_g),Id(L_h)]=0$,
where $Id(L_g)$ is the ideal of $L$ generated by $L_g$.
\end{lemma}

{\em Proof}. Any element of $I$ can be written
as a linear combination of products
$$
[x,y_1,\ldots,y_k,[z,t_1,\ldots,t_m]], \quad k\ge 0,m\ge 0,
$$
where $x\in L_g,z\in L_h$ and $y_1,\ldots,y_k,t_1,\ldots,t_m$
are homogeneous elements from $L$. Now the statement easily follows by Lemma~\ref{l1}.
\hfill $\Box$

As an immediate consequence of these remarks we get the following:

\begin{proposition}\label{p1}
Let $L=\oplus_{g\in G}L_g$ be a $G$-graded simple Lie algebra.
Then $Supp~L$ generates an abelian subgroup of
$G$.
\end{proposition}
\hfill $\Box$

In particular, if $L$ is simple in the
non-graded sense
then $Supp~L$ is a commutative subset of $G$.
Note that in general,
this property
does not hold  for semisimple algebras. For
example, let $L=B_1\oplus B_2$ be the direct sum of two
simple algebras isomorphic to $sl_2(F)$.
Given a group $G$ of order $2$, $G=\{e,g\}$,
we can define $G$-grading on $H=sl_2(F)$ by seting
$$
L_e=Span\left\{
\begin{pmatrix}
  1 &\; 0 \\
  0 &\; -1
\end{pmatrix}\right\},\quad
L_g=Span \left\{
\begin{pmatrix}
  0 &\; 1 \\
  0 &\; 0
\end{pmatrix},
\begin{pmatrix}
  0 &\; 0 \\
  1 &\; 0
\end{pmatrix}
\right\}.
$$
It is now sufficient to take any group $G$ with $gh\ne hg,~g,h\in G,~g^2=h^2=e$,
and define $G$-grading on $B_i,$
using $H_i$, $i=1,2$. Later we shall
show that all noncommutative gradings on semisimple Lie algebras are of
similar type.

\section{Structure of finite dimensional graded algebras}

Now let $G$ be a finite abelian group and let $F$ be an algebraically
closed field of characteristic zero. Recall a
well-known duality between
gradings and automorphism actions on $L$ (cf. e.g. \cite{BSZ2001}).
Let $\widehat G$ be the dual group for $G$, that is the group of all
irreducible characters on $G$. Since $G$ is finite abelian, the group
$\widehat G$ is isomorphic to $G$. If $L=\oplus_{g\in G}L_g$ is a
$G$-graded algebra then any $\chi\in\widehat G$ acts on $L$ by
the automorphism
$$
\chi*x_g=\chi(g)x_g,
$$
where $x_g\in L_g$ is a homogeneous element of degree $g$.

Clearly, all subspaces $L_g$ are stable under the
$\widehat G$-action.
Moreover, a subspace $V\subseteq L$ is graded if and only if $V$ is
$\widehat G$-stable. Conversely, if one defines the
$\widehat G$-action
on $L$ by automorphisms then $L$ can be decomposed into a direct sum
\begin{equation}\label{eq1}
L=\bigoplus_{g\in G}L_g
\end{equation}
where
$$
L_g=\{v\in L| ~\chi*v=\chi(g)v\quad{\rm for~all}\quad\chi\in \widehat G\}
$$
and the decomposition (\ref{eq1}) is a $G$-grading on $L$.

If $G$ is an infinite cyclic group generated by $t$ then we can also
define an action of the
infinite cyclic group generated by $\chi$ on $L=\oplus_{g\in G}L_g$,
by setting
$$
\chi^n* v=\lambda^{nk}v
$$
as soon as $v\in L_{t^k}$ where $\lambda\in F^*$ is a fixed element
of infinite order. As before, a subspace $V$ of $L$ is graded if and
only if $\chi*V=V$. If $L$ is finite dimensional and $Supp~L$ generates
$G$ then $G$ is a finitely generated abelian group and we can again
identify $G$-grading of $L$ with
the
$G$-action on $L$.

Using this duality we get the following result.

\begin{proposition}\label{p2}
Let $L=\oplus_{g\in G}$ be a finite dimensional $G$-graded
Lie algebra over an algebraically closed field of characteristic zero.
\begin{itemize}
\item[1)]
If $L$ is graded simple then $G$ is abelian and
$L=B_1\oplus\cdots\oplus B_n$ is semisimple
with isomorphic simple components $B_1,\ldots, B_n$.

\item[2)]
If $L$ is semisimple then $L=A_1\oplus\cdots\oplus A_m$ is a
direct sum of graded simple components and $Supp~A_i$ is a
commutative subset of $G,$ for any $1\le i\le m$.
\end{itemize}
\end{proposition}

{\em Proof}. First let $L$ be graded simple. Then $G$ is
abelian by Proposition~\ref{p1}. Suppose that
$L$ is not semisimple. Then due to the
duality between $G$-gradings and $G$-action its
solvable radical $R$ is a graded ideal since $R$
is stable under the action of any automorphism of
$L$. Hence $R=L$ and $L^2=[R,R]\ne L$ is also a
graded ideal and $L^2\ne 0$ by the definition of graded simplicity,
a contradiction. Hence $L$ is semisimple.
Consider the decomposition $L=B_1\oplus\cdots\oplus B_n$ into
a direct sum of simple ideals.

Let $f: L\rightarrow L$ be an automorphism of $L$. Clearly,
$f(B_1)=B_j$ for some $1\le j\le n$. Due to the duality between
$G$-grading and $G$-action, the orbit of $B_1$ under
the
$G$-action contains all summands $B_1,\ldots,B_n$ since $L$ is
graded simple. In particular, all $B_1,\ldots,B_n$ are isomorphic,
and we have proved assertion (1) of the proposition.

Now let $L$ be a semisimple $G$-graded algebra. As before,
$L$ is a direct sum of minimal ideals, $L=B_1\oplus\cdots\oplus B_n$.
Consider the
minimal graded ideal $A_1$ of $L$. If $A_1=L$ then there is
nothing to prove. If $A_1\ne L$ then $A_1$ is a
sum
of some $B_i$, say, $A_1=B_1\oplus\cdots\oplus B_t$. We shall
prove that $B_{t+1}\oplus\cdots\oplus B_n$ is also a
graded ideal of $L$.

First note that the centralizer of any homogeneous element
in $L$ is a graded subalgebra of $L$. Indeed, let
\begin{equation}\label{eq1a}
[a_h,b_{g_1}+\cdots+b_{g_m}]=0
\end{equation}
for some $a_h\in L_h, b_{g_1}\in L_{g_1},\ldots, b_{g_m}\in L_{g_m}$ with pairwise distinct $g_1,\ldots, g_m\in G$. If $[a_h,b_{g_i}]=0,$
for all $i=1,\ldots,m$ then we are done.
Otherwise fix all $g_{i_1},\ldots, g_{i_k}$ such that $c_1=[a_h,b_{g_{i_1}}]\ne 0,\ldots, c_k=[a_h,b_{g_{i_k}}]\ne 0$. Then $h$ commutes with $g_{i_1},\ldots,g_{i_k}$
by Lemma \ref{l1},
$c_1\in L_{hg_{i_1}},\ldots,c_k\in L_{hg_{i_k}}$ and all $hg_{i_1},\ldots,hg_{i_k}$
are distinct. On the other hand, $c_1+\cdots+c_k=0$, a contradiction. Hence the
equality (\ref{eq1a}) implies that all $b_{g_1},\ldots,b_{g_m}$ commute with $a_h$.

Now we observe
that the intersection of two graded subalgebras
is also a graded subalgebra. In particular, the centralizer of $A_1$ is a graded subapace
in $L$. However,
the centralizer of $A_1$ is $B_{t+1}\oplus\cdots\oplus B_n$
and we
complete
the proof of the proposition by  induction on dimension of $L$.
\hfill $\Box$

\begin{remark}\label{r1}
It follows by the previous proposition  that any finite dimensional
graded simple algebra is associated with one of the finite dimensional
simple Lie algebras.
In particular,
we can say that the graded simple algebra is an algebra one of the types $A_l~(l\ge 1),B_l~(l\ge 2),C_l~(l\ge 3),D_l~(l\ge 4), G_2,F_4, E_6,E_7,E_8$.
\end{remark}

Now we clarify the structure of a finite dimensional Lie algebra graded by a finite group.

\begin{proposition}\label{p3}
Let $G$ be a finite group and let $L=\oplus_{g\in G} L_g$ be a
finite dimensional $G$-graded Lie algebra over an algebraically
closed field of characeristic zero. Then the radical $R=Rad~L$ of
$L$ is a graded ideal and there exists a split extension $L=B+R$,
where $B$ is a maximal semisimple subalgebra of $L$ homogeneous in $G$-grading.
\end{proposition}

{\em Proof}. If $G$ is an abelian group then, as it was mentioned in the proof
of Proposition \ref{p2}, the solvable radical of $L$ is graded.

Now let $G$ be a noncommutative group. Consider a
decomposition $L=B+R,$ where $B$ is a maximal semisimple
subalgebra of $L$, $B=B_1\oplus\cdots\oplus B_n$ and
$B_1,\ldots, B_n$ are minimal ideals of $B$. By Proposition~\ref{p1},
algebra $L$ cannot be graded simple. Consider a maximal proper graded ideal $P$ of $L$.

First suppose that $B\subseteq P$. Then $\dim L/P=1$. Denote $R_0=Rad~ P$. Then
$R_0\subset R$ and $\dim R/R_0=1$. Moreover,
$R^2\subseteq R_0$ and $[B,R_0]\subseteq R_0$. In particular, $R_0$ is
an ideal of $L$. Since $\dim P<\dim L$,  we may suppose by  induction
that $R_0$ is a graded ideal of $P$ and $L$.
If $R_0\ne 0$ then also by induction, $R/R_0$ is a graded ideal of $L/R_0$,
hence $R$ is a graded ideal of $L$.

In case $R_0=0$ we have $\dim R=1$ and hence it is a
trivial $B$-module. Since $[R,R]=0$ we conclude that $R$
is the center of $L$. By Lemma \ref{l1},
the center of a graded
algebra is homogeneous and we
have completed the proof that $R$ is a graded ideal in this case.

Now let $B\not\subseteq P$. Then $L/P$ is graded simple and by Proposition~\ref{p2},
it is semisimple, i.e. $R\subseteq P$.
Applying the induction on dimension we again conclude that
$R$ is a homogeneous subspace of $P$ and hence of $L$.

Now we prove the existence of Levi subalgebra
homogeneous
in $G$-grading. In case $G$ is abelian we apply the duality
between $G$-grading and $G$-action by automorphisms on $L$.
By \cite{T1957} there exists a maximal semisimple subalgebra
$B$ stable under $\widehat G$-action. So, $B$ is a graded Levi
subalgebra of $L$.

In the
general case we consider graded factor-algebra $L/R$. By Proposition~\ref{p2} we have
$$
L/R=\bar A_1\oplus\cdots\oplus \bar A_m
$$
where any $\bar A_i$ is graded simple and $S_i=Supp~\bar A_i$ is a commutative subset of $G$.

Denote by $A_i$ the full preimage of $\bar A_i$ over $R$.
Then $A_i$ is graded and $Supp~A_i\supset S_i$.
Now we take a subalgebra $C_i$ of $A_i$ generated by all homogeneous
$x\in A_i$ with $\deg x\in S_i$. Since $S_i$ is commutative,
by the previous remark $C_i$ contains a homogeneous semisimple
subalgebra $B_i$ isomorphic $\bar A_i$. Note that $\dim(B_1+\cdots+B_m)=\dim L/R$ hence $B=B_1\oplus\cdots\oplus B_m$ is a homogeneous Levi subalgebra and we have
completed the proof.
\hfill $\Box$

Combining Propositions~\ref{p1}, \ref{p2} and \ref{p3} we immediately obtain the following:

\begin{theorem}\label{t1}
Let $L=\oplus_{g\in G} L_g$ be a finite dimensional Lie algebra over an
algebraically closed field of characteristic zero graded by a finite group $G$.
Then its solvable radical is homogeneous in $G$-grading and
there exists a Levi subalgebra $B$ homogeneous in $G$-grading. Moreover, $B$ is a direct sum $B=H_1\oplus\cdots\oplus H_m$ where any $H_i$ is graded
simple subalgebra with commutative support $Supp~H_i$ and in the
non-graded
sense $H_i$ is a
direct sum of isomorphic simple components.
\end{theorem}
\hfill $\Box$

\end{document}